\documentclass[10pt]{article}
\usepackage{pifont} % insert email package
\usepackage{graphicx,amssymb,mathrsfs,color,enumerate}
\date{}
\newfont{\bb}{msbm10}

\newcommand{\tr}{^{\sf T}}

\newtheorem{algorithm}{Algorithm}[section]
\newtheorem{remark}{Remark}[section]

\newtheorem{theorem}{Theorem}[section]

\newtheorem{lemma}{Lemma}[section]
\newtheorem{corollary}{Corollary}[section]

\def\x{\mathbf{x}}
\def\w{\mathbf{w}}
\def\X{\mathcal{X}}
\def\M{\mathcal{M}}
\def\R{\mathcal{R}}
\def\J{\mathcal{J}}

\begin{document}
\cleardoublepage \pagestyle{plain}
\bibliographystyle{plain}

\title{A family of multi-parameterized proximal point algorithms
\thanks{The work was  supported by  the Natural Science Foundation of China (11801455; 11571178) and the Fundamental Research Funds of China West Normal University (17E084; 18B031).}
   }

\author{
Jianchao Bai
\footnote{Department of Applied Mathematics,     Northwestern Polytechnical
University, Xi'an  710129,    China (\tt bjc1987@163.com).}
\quad   Ke Guo
\footnote{  (Corresponding author) School of Mathematics and Information, China West Normal University, Nanchong 637002,  China (\tt keguo2014@126.com).}
\quad   Xiaokai Chang
\footnote{  School of Science, Lanzhou University of Technology, Lanzhou 730050,   China (\tt xkchang@lut.cn).}
}
\maketitle

\centerline{\small\it\bf Abstract}\vskip 1mm
In this paper, a multi-parameterized proximal point algorithm combining with a relaxation step is developed   for solving  convex minimization problem subject to linear constraints. We show its    global convergence and sublinear convergence rate from the prospective of  variational inequality.   Preliminary numerical experiments on testing a sparse minimization problem   from signal processing indicate that the proposed algorithm performs better than some well-established methods.


% 65Y20---Complexity and performance of numerical algorithms
% 90C25 Convex programming
% 92C55---Biomedical imaging and signal processing
\bigskip

%====================================================================%
\section{Introduction}
%====================================================================%
We focus on the following convex minimization problem with linear equality constraints,
\begin{equation} \label{Sec1-Prob}
\begin{array}{lll}
\min \{  f(\x)\ |   A\x=b, \x\in \X\},
\end{array}
\end{equation}
where $f:\R^{n}\rightarrow\R\cup\{+\infty\}$ is a  proper closed convex function  but possibly nonsmooth, $A\in\R^{m\times n}$ and $ b\in\R^{m}$ are given  matrix and vector, respectively,
$\X\subseteq \R^{n}$ is a    closed convex set.
Without loss of generality,
the solution set  of the problem (\ref{Sec1-Prob}) denoted by $\X^*$   is assumed to be nonempty.

The augmented Lagrangian method (ALM), independently proposed by Hestenes  \cite{H1969} and Powell \cite{P69}, is a benchmark method  for solving problem (\ref{Sec1-Prob}). Its iteration scheme reads as
 \[
 \left \{\begin{array}{l}
\x^{k+1}=\arg\min\limits_{\x\in\X} f(\x)+\frac{\beta}{2}\|A\x-b-\lambda^k/\beta\|^2,\\
\lambda^{k+1}=\lambda^k- \beta(A\x^{k+1}-b),
\end{array}\right.
\]
where  $\beta,\lambda$ denote  the penalty parameter and  the Lagrange multiplier w.r.t. the equality constraint, respectively.  As analyzed in \cite{Roc76},    ALM  can be viewed as an
 application of  the well-known proximal point algorithm (PPA) that can date back to the
 seminal work of Martinet \cite{M70} and Rockafellar \cite{R76} for  the dual problem of  (\ref{Sec1-Prob}).
 Obviously, the efficiency of   ALM heavily depends on the solvability of the $\x$-subproblem, that  is,  whether or not the core $\x$-subproblem has closed-form solution.  Unfortunately, in many real applications \cite{bAIZ18,DonoT08,KKL07,MaNi2016},  the coefficient matrix $A$ is not  identity matrix (or does not satisfies $AA\tr=I_m$), which makes it difficult even infeasible for solving this  subproblem of ALM.  To overcome such  difficulty, Yang and Yuan \cite{YY13} proposed a  linearized ALM aiming at linearizing the $\x$-subproblem such that  its closed-form solution    can be easily derived.  We refer to the recent progress on this direction \cite{bjkj9,HMY19}.

 Under basic regularity condition  $\mbox{ri~dom}(f)\cap X\neq\emptyset$, it is well-known that $x^{*}$ is an optimal solution of  (\ref{Sec1-Prob}) if and only if there exists $\lambda^{*}\in \R^{m}$
 such that  the following  variational inequality   holds
\begin{equation} \label{Sec2-003}
\textrm{VI}(f, \J,\M): \quad  f(\x)-f(\x^*)+ (w-w^*)\tr \J(w^*)\geq 0,\quad\forall w\in\M,
\end{equation}
where
\[
w=\left(\begin{array}{c}
\x\\  \lambda\\
\end{array}\right),\quad
 w^*=\left(\begin{array}{c}
\x^*\\  \lambda^*\\
\end{array}\right),\quad
\J(w)=\left(\begin{array}{c}
-A\tr\lambda\\
 A\x-b
\end{array}\right)\quad \textrm{and}\quad
\M=\X\times \R^{m}.
\]

From the aforementioned    assumption on the   solution set of the problem (\ref{Sec1-Prob}),
 the solution set of (\ref{Sec2-003}) denoted by $\M^*$  is also nonempty.
 When  PPA is applied to solve the variational inequality $\textrm{VI}(f, \J,\M)$, it usually reads  the unified updating scheme: at the $k$-th iteration, find $w^{k+1}$ satisfying
 \begin{equation} \label{Sec2-3}
f(\x)-f(\x^{k+1})+ (w-w^{k+1})\tr \left[\J(w^{k+1})+G(w^{k+1}-w^k)\right]\geq 0,\quad\forall w\in\M.
\end{equation}
We call $G$  the proximal matrix  that is usually required to be symmetric  positive definite to  ensure  the convergence of (\ref{Sec2-3}). To our knowledge, this idea was initialized by He et al. \cite{HLHY02}.
 Clearly, different structures of $G$ would result in different versions of PPA.  From a computational perspective,  our motivation is to design a multi-parameterized PPA  for solving problem (\ref{Sec1-Prob}) while maintaining the efficiency as the linearized ALM, although the   feasible starting point may be  different. Interestingly, many customized proximal matrices  shown in \cite{GuHeYuan2014,HeYuanZhang2013,MaNi2016,ZWY15} turn out to be special cases of our multi-parameterized proximal matrix (See Remark \ref{agaga} for details).  In this sense, our  proposed  algorithm can be viewed as a  general customized PPA for solving problem (\ref{Sec1-Prob}).
 Moreover, we adopt a relaxation strategy to accelerate the convergence.

\section{Main Algorithm}

In this paper, we  design the following  multi-parameterized proximal matrix
\begin{equation} \label{Sec2-G}
G=\left[\begin{array}{cc}
rI_n +\frac{(\theta-1)^2-\rho}{s}A\tr A &   (\theta-1) A\tr\\
(\theta-1) A &  sI_m
\end{array}\right]\in \R^{(n+ m)\times (n+ m)},
\end{equation}
where   $I_m\in \R^{m\times m}$   denotes the identity matrix, $\theta$ is an arbitrary real scalar and
\begin{equation} \label{Sec2-region}
\rho\in(-\infty,1], \quad (r,s)\in \left\{(r,s)|~ r>0, s>0, rs>\|A\tr A\|_2\right\}.
\end{equation}
The notation   $\|A\tr A\|_2=\sqrt{\lambda_{\max}(A\tr A)}$ represents the  spectral norm of $A\tr A.$  It is easy to check that the above matrix $G$ is symmetric positive definite for any  parameters  $(\rho,r,s)$ satisfying (\ref{Sec2-region}).

Now,   substituting the   matrix $G$ into (\ref{Sec2-3})  we have
 \begin{equation} \label{Sec2-4}
 \left \{\begin{array}{l}
A\x^{k+1}-b + (\theta-1) A(\x^{k+1}-\x^k) +s(\lambda^{k+1}-\lambda^k)=0,\\
\x^{k+1}\in\X, \quad f(\x)- f(\x^{k+1}) + (\x-\x^{k+1})\tr R_{k+1}\geq 0,\quad\forall \x\in\X,
\end{array}\right.
\end{equation}
with
\begin{equation} \label{Sec2-5}
R_{k+1}= -A\tr\lambda^{k+1} + \left[rI_n +\frac{(\theta-1)^2-\rho}{s}A\tr A\right](\x^{k+1}-\x^k)+  (\theta-1)A\tr(\lambda^{k+1}-\lambda^k).
\end{equation}
By the equation in  (\ref{Sec2-4}), it can be deduced  that
\[
\lambda^{k+1}=\lambda^k-\frac{1}{s}\left[ \theta(A\x^{k+1}-b) +(1-\theta)(A\x^k-b)\right],
\]
which further makes (\ref{Sec2-5}) become
\begin{eqnarray*}
R_{k+1}&=& -A\tr\left[(2-\theta)\lambda^{k+1} +(\theta-1)\lambda^k+ \frac{\rho-(\theta-1)^2}{s} A(\x^{k+1}-\x^k)\right]+r(\x^{k+1}-\x^k)\\
% &=&-(2-\theta)A\tr\left[\lambda^k-\frac{1}{s}\left[ \theta(A\x^{k+1}-b) +(1-\theta)(A\x^k-b)\right]\right] +r(\x^{k+1}-\x^k)\\
%&&-A\tr\left[(\theta-1)\lambda^k+ \frac{\rho-(\theta-1)^2}{s} A(\x^{k+1}-\x^k)\right]\\
&=& -A\tr\left[\lambda^k- \frac{(2-\theta)}{s}(A\x^k-b)\right]+\left[rI_n+\frac{\rho-1}{s}A\tr A\right](\x^{k+1}-\x^k).
\end{eqnarray*}
Based on   the inequality  in   (\ref{Sec2-4}),  i.e., the first-order optimality condition of $\x$-subproblem, we obtain
\begin{equation}\label{key}
\x^{k+1}=\arg\min\limits_{\x\in\X}\left\{
f(\x) +\frac{r}{2}\left\|\x-\x^k -\frac{1}{r}A\tr\left[\lambda^k-\frac{2-\theta}{s}(A\x^k-b)\right]\right\|^2+
\frac{\rho-1}{2s}\left\|A(\x-\x^k)\right\|^2
\right\}.
\end{equation}
Then,    our  relaxed  multi-parameterized PPA (RM-PPA) is described as Algorithm \ref{algo1}, where we use $\widetilde{w}^{k}:=(\widetilde{\x}^k,\widetilde{\lambda}^k)$ to replace the  output of  (\ref{Sec2-3}) with given iterate $(\x^k,\lambda^k)$,  and we use $(\x^{k+1},\lambda^{k+1})$ to stand for the new iterate  after combining a relaxation step. Finally, the inequality (\ref{Sec2-3}) becomes
 \begin{equation} \label{Sec2-31}
f(\x)-f(\widetilde{\x}^k)+ (w-\widetilde{w}^k)\tr \left[\J(\widetilde{w}^k)+G(\widetilde{w}^k-w^k)\right]\geq 0,\quad\forall w\in\M.
\end{equation}
\begin{algorithm}\label{algo1}
\vskip1mm
\hrule\vskip1.5mm
[RM-PPA for solving problem (\ref{Sec1-Prob})]\vskip0.5mm
\noindent \verb"1  "\verb"Choose" $\sigma\in(0,2), \theta\in \R$ \verb"and" $ (\rho,r,s)$ \verb"satisfying" (\ref{Sec2-region}).\\
\verb"2  Initialize" $(\x^0,\lambda^0)\in \M.$\\
\verb"3  for" $k=1,2,\cdots,$ \verb"do" \\
\verb"4"\indent\quad\ \
$\widetilde{\x}^{k}=\arg\min\limits_{\x\in\X}\left\{
f(\x) +\frac{r}{2}\left\|\x-\x^k -\frac{1}{r}A\tr\left[\lambda^k-\frac{2-\theta}{s}(A\x^k-b)\right]\right\|^2+
\frac{\rho-1}{2s}\left\|A(\x-\x^k)\right\|^2
\right\}$.\\
\verb"5" \indent\quad\
$\widetilde{\lambda}^{k}=\lambda^k-\frac{1}{s}\left[ \theta(A\widetilde{\x}^k-b) +(1-\theta)(A\x^k-b)\right]$.\\
\verb"6" \indent\quad\
$\left(\begin{array}{c}
\x^{k+1}\\  \lambda^{k+1}\\
\end{array}\right)=
\left(\begin{array}{c}
\x^k\\  \lambda^k\\
\end{array}\right)-\sigma
\left(\begin{array}{c}
\x^k-\widetilde{\x}^k\\  \lambda^k-\widetilde{\lambda}^k\\
\end{array}\right)$.\\
\verb"7  end"
\vskip2mm\hrule\vskip3mm
\end{algorithm}

\begin{remark} If we set $\rho=1$ in (\ref{key}),  then  the $\x$-subproblem  amounts to estimating the proximity operator of $f$ when $X=\R^{n}$. The implementation of (\ref{key}) for such cases is thus extremely simple. Here, we allow $\rho\in (-\infty,1]$  just from   the theoretical point of view.

\end{remark}

\begin{remark}  \label{agaga}
Note that  ${1}/{s}$ in step 5 actually plays a role of penalty parameter in ALM, while $r$ can be treated as the proximal parameter as used in the customized PPA \cite{HeYuanZhang2013}. The quadratic term \[\frac{\rho-1}{2s}\left\|A(\x-\x^k)\right\|^2=\frac{\rho-1}{2s}\left\|(A\x-b)-(A\x^k-b)\right\|^2\]
plays a second penalty role for the equality constraint relating to its   $k$-th iteration.
By the way of updating  $\widetilde{\lambda}^{k}$, it uses the convex combination of the feasibility error at the current iteration and the former iteration when $\theta\in[0,1]$. The parameterized   matrix designed in this paper is more general than some in the literature:
 \begin{itemize}
 \item
If $(\theta,\rho)=(0,1)$, then our   matrix $G$ given by (\ref{Sec2-G}) will become  that in \cite[Eq.(2.5)]{HeYuanZhang2013}. And in such case, the  variable $\lambda$ updates practically in the same way as in ALM, but the core $\x$-subproblem in Algorithm \ref{algo1}  is a proximal mapping to have a unique  minimum    since the subproblem  is strongly convex.
\item
If $(\theta,\rho)=(2,1)$, then our parameterized proximal matrix   turns to  the matrix $Q$ involved in \cite[page 158]{GuHeYuan2014}. If $(\theta,\rho)=(\tau+1,1)$,  then our   matrix $G$   is  identical to that in \cite[Eq. (3.1)]{MaNi2016} but Algorithm \ref{algo1} uses an additional relaxation step for fast convergence. Moreover, we establish the worst-case  $\mathcal{O}(1/t)$ ergodic convergence rate for the  objective function value  error  and the feasibility error.
\item
 Regardless of step 6, it is easy to check that   Algorithm \ref{algo1}  with $\theta=\rho=1$  is a linearization of ALM:
\begin{equation} \label{linear-ALM}
 \left \{\begin{array}{l}
\x^{k+1}=\arg\min\limits_{\x\in\X} \left\{f(\x)+\frac{1}{2s}\|A\x-b-s\lambda^k\|^2+\frac{1}{2}\|\x-\x^k\|^2_{rI_n-\frac{1}{s}A\tr A}\right\},\\
\lambda^{k+1}=\lambda^k- \frac{1}{s}(A\x^{k+1}-b).
\end{array}\right.
\end{equation}
Specifically, by letting $\beta=1/s$  the  scheme (\ref{linear-ALM}) is  ALM with extra   proximal term $\frac{1}{2}\|\x-\x^k\|^2_{rI_n-\frac{1}{s}A\tr A}$ which eliminates  the   term $\|A\x\|^2 $ in the iteration. Algorithm \ref{algo1}, in such choice of parameters,    is   a   linearized  ALM.  Besides, our   parameter $\theta$ is more general and flexible than that $\theta\in[-1,1]$    in \cite{ZWY15}.
\end{itemize}
\end{remark}

\section{Convergence   Analysis}
Before  analyzing  the global convergence and sublinear convergence rate of Algorithm \ref{algo1}, we give a fundamental lemma as the following.

\begin{lemma}\label{key-lem}
The sequence  $\{w^k\}$ generated by Algorithm \ref{algo1} satisfies
\[
\|w^{k+1}-w^*\|_{\widetilde{G}}^2+ \frac{1-\rho}{s}\|A\x^{k+1}-A\x^*\|^2\leq \|w^k-w^*\|_{\widetilde{G}}^2+\frac{1-\rho}{s}\|A\x^k-A\x^*\|^2 -T_{k+1}, ~~\forall w^*\in\M^*,
\]
where  $T_{k+1}$  is given by (\ref{S2-04}) and
\begin{equation} \label{Sec2-barG}
\widetilde{G}=\left[\begin{array}{cc}
rI_n +\frac{(\theta-1)^2-1}{s}A\tr A &   (\theta-1) A\tr\\
(\theta-1) A &  sI_m
\end{array}\right].
\end{equation}
\end{lemma}
\noindent{\bf Proof } According to  the inequality (\ref{Sec2-003}) and the   skew-symmetric property of  $\J(w)$, i.e.
\begin{equation} \label{S2-0012}
(w-\bar{w})\tr[\J(w)-\J(\bar{w})]\equiv0,\quad \textrm{for any } w, \bar{w}\in\M,
\end{equation}
  the inequality (\ref{Sec2-31}) with setting $w=w^*$  gives
$(w^*-\widetilde{w}^k)\tr G(\widetilde{w}^k-w^k)\geq 0.$
Note that the step 6 shows
\begin{equation} \label{S2-00}
\widetilde{w}^k-w^k=(w^{k+1}-w^k)/\sigma,
\end{equation}
so we have
\[
(w^*-\widetilde{w}^k)\tr G(w^{k+1}-w^k)\geq 0.
\]
Since the matrix $G$ can be decomposed as
$
G=\widetilde{G} +\textrm{Diag}\left((1-\rho)A\tr A/s, \textbf{0}_m\right),
$
where  $\textbf{0}_m$ denotes the  zero matrix of size $m\times m$  and $\widetilde{G}$ is given by (\ref{Sec2-barG}),   we thus obtain
\begin{equation} \label{S2-01}
(w^*-\widetilde{w}^k)\tr\widetilde{G}( w^{k+1}-w^k )+\frac{1-\rho}{s}( A\x^*- A\widetilde{\x}^k)\tr( A\x^{k+1}-A\x^k)\geq 0.
\end{equation}
Then,   applying the identity \[2(a-l)\widetilde{G}(c-d)= \|a-d\|_{\widetilde{G}}^2-\|a-c\|_{\widetilde{G}}^2 + \|c-l\|_{\widetilde{G}}^2-\|d-l\|_{\widetilde{G}}^2 \]
 to the left-hand side of  (\ref{S2-01}), the following inequality holds immediately
\begin{equation} \label{S2-02}
 \|w^*-w^{k+1}\|_{\widetilde{G}}^2+\frac{1-\rho}{s}\| A\x^*-A\x^{k+1}\|^2\leq  \|w^*-w^k\|_{\widetilde{G}}^2+\frac{1-\rho}{s}\| A\x^*-A\x^k\|^2- T_{k+1},
\end{equation}
 where
\[
T_{k+1}=  \|w^k-\widetilde{w}^k\|_{\widetilde{G}}^2 -\|w^{k+1}-\widetilde{w}^k\|_{\widetilde{G}}^2
+\frac{1-\rho}{s}\left( \|A(\x^k-\widetilde{\x}^k)\|^2- \|A(\x^{k+1}-\widetilde{\x}^k)\|^2 \right).
\]
Substituting (\ref{S2-00}) into the expression of $T_{k+1}$, it can be deduced that
\begin{eqnarray} \label{S2-04}
T_{k+1}&=&  \|w^k-\widetilde{w}^k\|_{\widetilde{G}}^2 -\|w^{k+1}-w^k+w^k-\widetilde{w}^k\|_{\widetilde{G}}^2
+\frac{1-\rho}{s}\left( \|A(\x^k-\widetilde{\x}^k)\|^2- \|A(\x^{k+1}-\x^k+\x^k-\widetilde{\x}^k)\|^2 \right)\nonumber\\
&=&\sigma(2-\sigma)\|w^k-\widetilde{w}^k\|_{\widetilde{G}}^2
+\frac{(1-\rho)\sigma(2-\sigma)}{s}\|A(\x^k-\widetilde{\x}^k)\|^2 \nonumber\\
&=&\frac{2-\sigma}{\sigma}\|w^k-w^{k+1}\|_{\widetilde{G}}^2
+\frac{(1-\rho)(2-\sigma)}{s\sigma}\|A(\x^k-\x^{k+1})\|^2.
\end{eqnarray}
This completes the whole proof. $\ \ \ \blacksquare$

Lemma \ref{key-lem} shows the sequence $\{w^k\}$ is contractive under the $\widetilde{G}$-norm w.r.t. the solution set $\M^*$, since the matrix $\widetilde{G}$ is positive definite and the term $T_{k+1}\geq0$. Similar to the convergence proof in e.g. \cite{bAIZ18} and the proof of Lemma \ref{key-lem}, the  global convergence and sublinear convergence rate of Algorithm \ref{algo1}   can be easily established as as below, whose proof is omitted here for the sake of conciseness.

\begin{theorem} \label{Sec3-conver}
Let $ (\rho,r,s)$ satisfy (\ref{Sec2-region}) and
 $\{w_k\}$ be generated by Algorithm \ref{algo1}. Then,
 \begin{itemize}
 \item
  there exists a $w^\infty \in \mathcal{M}^*$ such that $\lim_{k \to \infty} w^k = w^\infty$;
\item
  for any $t>0$, let $\w_t=\frac{1}{t+1}\sum_{k=0}^{t}\widetilde{w}^k$ and $\x_t=\frac{1}{t+1}\sum_{k=0}^{t}\widetilde{\x}^k$. Then,
  \[
  f(\x_t)-f(\x)+ (\w_t-w)\tr \J(w)\leq \frac{1}{2\sigma(t+1)}\left\{\|w^0-w\|_{\widetilde{G}}^2+ \frac{1-\rho}{s}\|A\x^0-A\x\|^2\right\},\quad\forall w\in\M.
  \]
\end{itemize}
\end{theorem}

Theorem \ref{Sec3-conver} illustrates that Algorithm \ref{algo1} converges globally with a sublinear ergodic convergence rate. Furthermore,  we can deduce a compact result   as the following corollary by making using of the second result in Theorem \ref{Sec3-conver}. For any $\xi>0$, let $\Gamma_\xi=\{\lambda~ | ~ \xi\geq\|\lambda\|\}$ and
\begin{equation} \label{Sec2-b11}
\gamma_\xi=\inf\limits_{\x^*\in\X^*}\sup\limits_{\lambda\in \Gamma_\xi}\left\|(\x^0-\x^*;\lambda^0-\lambda)\right\|_{\widetilde{G}}^2+\frac{1-\rho}{s}\|A\x^0-b\|^2.
\end{equation}

\begin{corollary} \label{Sec3-coll}
Let
 $\{w_k\}$ be generated by Algorithm \ref{algo1}. For any $\xi>0$, there exists a $\gamma_\xi<\infty$ defined in (\ref{Sec2-b11}) such that for any $t>0$, we have
 \[
 f(\x_t)-f(\x^*)+ \xi \|A\x_t-b\|\leq \frac{\gamma_\xi}{2\sigma(t+1)}, \quad\forall \x^*\in\X^*.
 \]
\end{corollary}
\noindent{\bf Proof }  By  making use of   the  identity in (\ref{S2-0012})
and by setting $w=(\x^*,\lambda)\in\X^*\times \R^m$ into the second result of Theorem \ref{Sec3-conver}, we have
\begin{eqnarray}\label{Sec2-b01}
&&f(\x_t)- f(\x^*)+ (\w_t-w)\tr \J(w)\nonumber\\
&=&f(\x_t)- f(\x^*)- \lambda\tr A(\x_t-\x^*)+ (\lambda_t-\lambda)\tr (A\x^*-b)\nonumber\\
&=&f(\x_t)- f(\x^*)- \lambda\tr (A\x_t-b)\leq   \frac{1}{2\sigma(t+1)}\left\{\left\|(\x^0-\x^*;\lambda^0-\lambda)\right\|_{\widetilde{G}}^2+ \frac{1-\rho}{s}\|A\x^0-b\|^2\right\},
\end{eqnarray}
where the second equality and the final inequality use  $A\x^*=b$. Then, it follows from (\ref{Sec2-b01}) that
\begin{eqnarray*}
f(\x_t)- f(\x^*)+ \xi\| A\x_t-b\|&=& \sup\limits_{\lambda\in \Gamma_\xi}[f(\x_t)- f(\x^*)- \lambda\tr (A\x_t-b)]
\\
&\leq&  \frac{1}{2\sigma(t+1)}\left\{\inf\limits_{\x^*\in\X^*}\sup\limits_{\lambda\in \Gamma_\xi}\left\|(\x^0-\x^*;\lambda^0-\lambda)\right\|_{\widetilde{G}}^2+\frac{1-\rho}{s}\|A\x^0-b\|^2\right\},
\end{eqnarray*}
which, by the definition of $\gamma_\xi$ in (\ref{Sec2-b11}), completes the proof. $\ \ \ \blacksquare$

In a similar  analysis  to (\ref{Sec2-b01}) together with (\ref{Sec2-003}),  we can  derive  $f(\x_t)- f(\x^*)- (\lambda^*)\tr (A\x_t-b)\geq 0$   showing that $f(\x_t)- f(\x^*)  \geq - \|\lambda^*\| \|A \x_t - b\|.$ So,
taking $\xi = 2 \|\lambda^*\|+1$ in Corollary \ref{Sec3-coll}, the following inequality
\[
(\|\lambda^*\|+1)  \|A \x_t - b\| \leq  f(\x_t)- f(\x^*) + (2 \|\lambda^*\|+1)  \|A \x_t -b \|
\le \frac{\bar{\gamma}}{2\sigma(t+1)},
\]
holds with $\bar{\gamma} = \gamma_\xi < \infty$   given by (\ref{Sec2-b11}). Rearranging the above inequality, we have
\begin{equation} \label{xyz-fea-error}
\|A \x_t - b\| \leq   \frac{1}{ 2\sigma(t+1)}\ \frac{ \bar{\gamma}}{\|\lambda^*\|+1}.
\end{equation}
Hence, we will also have $ f(\x_t)- f(\x^*) \geq - \|\lambda^*\| \|A \x_t - b\| \ge
- \frac{ \bar{\gamma}}{2\sigma (t+1)}$ showing that
\begin{equation} \label{xyz-f-error}
| f(\x_t)- f(\x^*)| \le \frac{ \bar{\gamma}}{2 \sigma(t+1)}.
\end{equation}
According to  (\ref{xyz-f-error}) and  (\ref{xyz-fea-error}), both the objective function value error and the feasibility error at the
ergodic iterate $\x_t$ will decrease in the order of $\mathcal{O}(1/t)$ as $t$ goes to infinity.

%============================================================================
\section{Numerical Experiments}\label{numex}
%============================================================================
In this section,  we  apply   the proposed algorithm   to  solve the following  $l_1$-minimization problem  from signal processing \cite{KKL07}, which aims   to reconstruct a  length $n$ sparse signal from $m(<n)$ observations:
\begin{equation} \label{Sec4-examp}
\begin{array}{lll}
\min \left\{  \|\x\|_1\ |   A\x=b, \x\in \R^n\right\}.
\end{array}
\end{equation}
Note that this is a special case of (\ref{Sec1-Prob}) with specifications $f=\|\x\|_1$ and $\mathcal{X}=\R^n$.  Applying Algorithm \ref{algo1}  to problem (\ref{Sec4-examp}), we derive\footnote{The proximity operator is defined as $\textrm{Prox}_{f,r}(\x)=\arg\min \left\{f(\x)+\frac{r}{2}\|\x-c\|^2~| c\in\R^n\right\}.$}     $\x^{k+1}=\textrm{Prox}_{\|\x\|_1,r}(\x)$
that can be explicitly expressed by the   shrinkage operator \cite{DonoT08} to  be coded by the MATLAB inner function `\verb"withresh"'.   Followed by Lemma \ref{key-lem},     we use the following  stopping criterions   under given   tolerance:
\begin{equation} \label{stop}
\textrm{It\_{err}(k)}:= \frac{\max\left\{\|\x_k-\x_{k-1}\|, \|\lambda_k-\lambda_{k-1}\|\right\}}{\max\left\{\|\x_{k-1}\|, \|\lambda_{k-1}\|,1\right\}}\leq 10^{-4} \quad \textrm{and}  \quad \textrm{Eq\_err(k)}:=\frac{\|A\x^k-b\|}{\|b\|}\leq 10^{-4}.
\end{equation}
All of the forthcoming experiments use the same starting points $(\x^0,\lambda^0)=(0,0)$ and  are tested   in MATLAB R2018a (64-bit) on Windows 10 system  with an Intel Core i7-8700K CPU (3.70 GHz)  and 16GB memory.

Consider   an original signal $\x\in \mathcal{R}^{10000}$   containing 180 spikes with amplitude $\pm 1$. The measurement matrix $A\in \mathcal{R}^{3000\times 10000}$ is drawn  firstly from  the standard
norm distribution $\mathcal{N}(0,1)$ and then each of its row is normalized. The observation $b$ is generated by $b=A\x+v$, where $v$ is generated by the   Gaussian distribution $\mathcal{N}(0,0.01^2I)$ on $\mathcal{R}^{3000}$. With  the tuned parameters $(r,s,\rho,\sigma)=(8, 1.01\|A\tr A\|_2/r,1,1.4)$, some   computational results under different parameter $\theta$ are shown in Table 1  in which we present the  number of iterations (IT), the CPU time in seconds (CPU), the final obtained residuals $\textrm{It\_{err}}$ and $\textrm{Eq\_err}$, as well as the recovery error $\textrm{RE}={\|\x^k-\x_{\textrm{orig}}\|}/{\|\x_{\textrm{orig}}\|}$. Reported results from Table 1 indicate that the choice of   $\theta$ could make a great effect on the performance of our algorithm w.r.t. IT and CPU. And it seems that setting $\theta=0.5$ would be a reasonable choice to save the CPU time and to cost fewer  number of iterations.  The reconstruction results under $\theta=0.5$ are shown in Fig. 1, from which    the solution obtained by our algorithm always has the correct number of pieces and is  closer to the original noseless signal.

\vskip5mm\begin{center}
\begin{tabular}{cccccc}
\hline
 $\theta$ &  IT  & CPU    &  \textrm{It\_{err}} & \textrm{Eq\_{err}} & \textrm{RE}\\
\hline
-5& 886	  & 37.40  & 9.97e-5    & 7.96e-5     &  6.93e-2    \\
-2& 827 &  34.85  & 9.99e-5    &  8.34e-5    &   6.92e-2     \\
-1& 844 &  34.76  & 9.96e-5    &  8.52e-5    &   6.92e-2     \\
-0.5&851&  34.83  &   9.98e-5  &  8.61e-5    & 6.92e-2      \\
0 &  845 &  34.50  &  9.97e-5  & 8.61e-5    &      6.92e-2 \\
0.2 &851&   35.15  &  9.99e-5   & 8.62e-5     &    6.92e-2   \\
0.5&  826&  33.93  &  9.98e-5   & 8.62e-5     &     6.91e-2  \\
1  &  840&   34.64 &  9.97e-5   & 8.65e-5    &     6.91e-2 \\
2  &  832&  34.15 & 9.99e-5    & 8.65e-5    &    6.91e-2  \\
5  &  855&  35.68 & 9.99e-5    & 8.43e-5    &     6.91e-2 \\
10 & 881	 &  36.63  & 9.93e-5    & 7.84e-5    &   6.92e-2  \\
\hline
\end{tabular}
\end{center}
\begin{center}
Table 1:\ Results by Algorithm \ref{algo1} with different parameter $\theta$.
\end{center}

\begin{figure}[htbp]
 \begin{minipage}{1\textwidth}
 \def\figurename{\footnotesize Fig.}
 \centering
\resizebox{17cm}{3cm}{\includegraphics{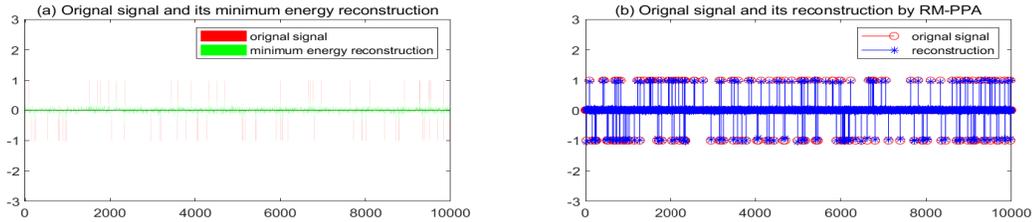} }
\caption{\footnotesize Comparison results by its minimum energy reconstruction (a) and by  RM-PPA  with  $\theta=0.5$ (b).}
   \end{minipage}
\end{figure}

\begin{figure}[htbp]
 \begin{minipage}{1\textwidth}
 \def\figurename{\footnotesize Fig.}
 \centering
\resizebox{17cm}{5cm}{\includegraphics{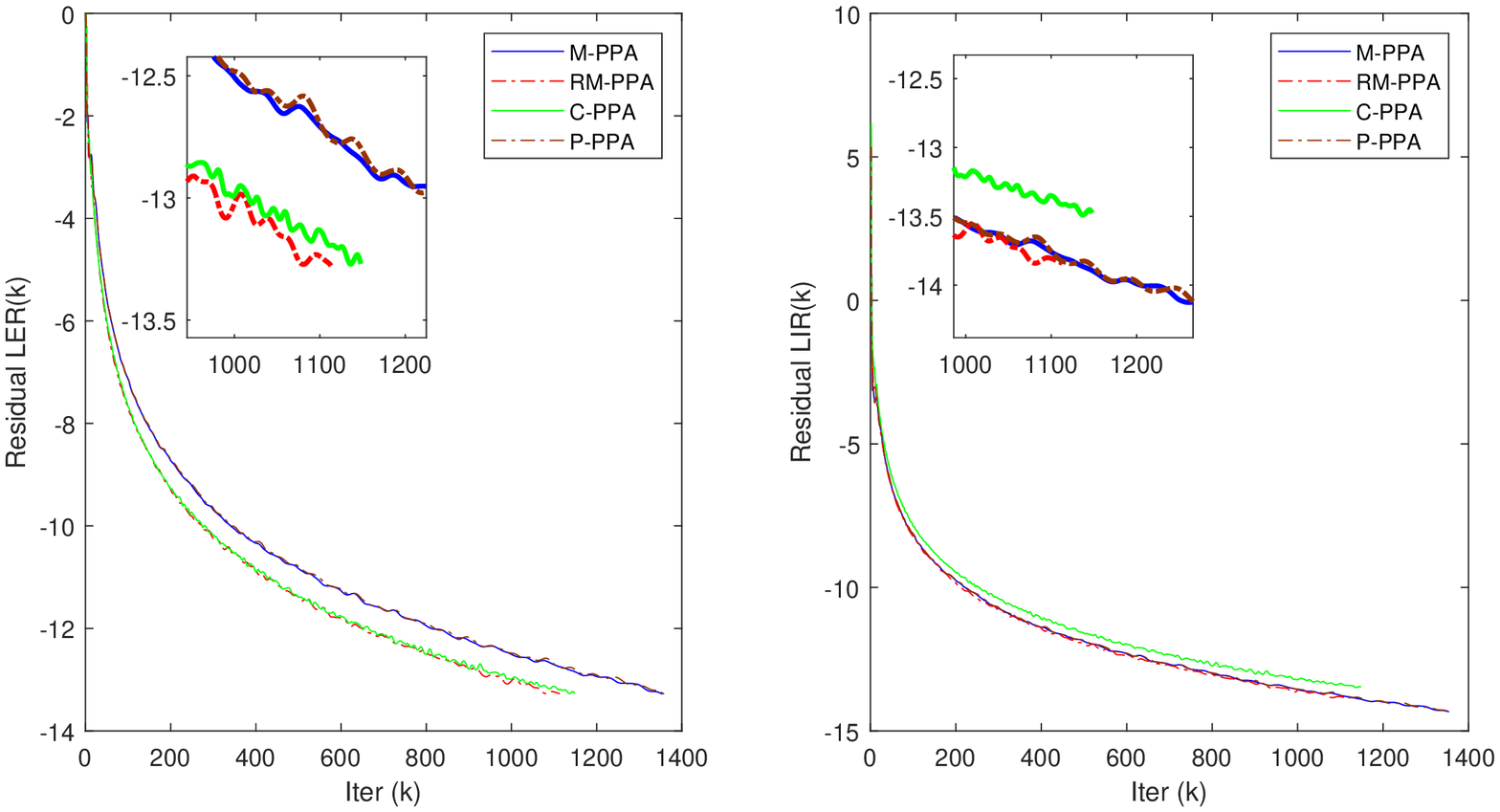} }
\caption{\footnotesize Convergence behaviors   of the residuals LER and LIR by different algorithms.}
   \end{minipage}
\end{figure}

Next, we consider the problem (\ref{Sec4-examp}) with large-scale dimensions $(m,n)=(3000, 20000)$. By comparing the  proposed Algorithms \ref{algo1} (RM-PPA) with the aforementioned tuned parameters to
 \begin{itemize}
  \item
  Algorithm \ref{algo1} without the relaxation step (``M-PPA"),
 \item
The customized PPA (``C-PPA",  \cite{HeYuanZhang2013}) with parameters $(\gamma,r,s)=(1.8, 8, 1.02\|A\tr A\|_2/r)$,
 \item
 The parameterized PPA (``P-PPA",  \cite{MaNi2016}) with parameters $(t,r,s)=(-1, 8, 1.02\|A\tr A\|_2/r)$,
\end{itemize}
 we show   comparative results about  the convergence behaviors  of the residuals $\textrm{LER(k)}:=\log_2(\textrm{Eq\_{err}(k)})$ and $\textrm{LIR(k)}:=\log_2(\textrm{It\_{err}(k)})$ in Fig. 2, respectively. The effect on recovering the original signal with different algorithms is shown in Fig. 3. Here, we   emphasize that the  parameter values in \cite{HeYuanZhang2013,MaNi2016} can not terminate the algorithms  C-PPA and P-PPA  because of the fact $\|A\tr A\|_2=1$ for their examples, so we set $r$ the same value as ours but keep $s$ as the value  in their experiments. From Figs. 2-3, we observe that M-PPA is competitive to P-PPA and RM-PPA (that is, Algorithm \ref{algo1}) performs better than the rest three algorithms.

\begin{figure}[htbp]
 \begin{minipage}{1\textwidth}
 \def\figurename{\footnotesize Fig.}
 \centering
\resizebox{17cm}{6cm}{\includegraphics{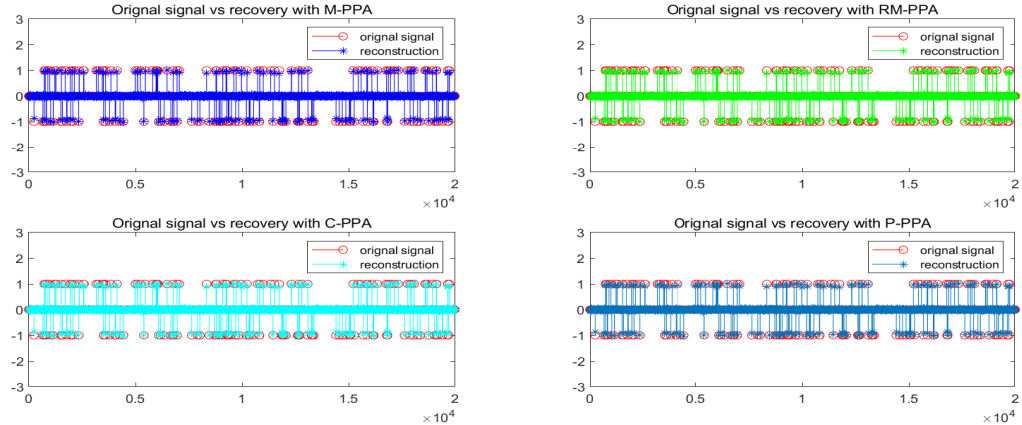} }
\caption{\footnotesize Comparison results of the problem (\ref{Sec4-examp}) with $n=20000$ by different algorithms.}
   \end{minipage}
\end{figure}

\end{document}